\documentclass[preprint,12pt]{elsarticle}
\usepackage{amsthm,amsmath,amssymb}
\usepackage[colorlinks=true,citecolor=black,linkcolor=black,urlcolor=blue]{hyperref}
\usepackage{graphicx}
\usepackage{float}
\usepackage{mathtools}
\usepackage{epstopdf}
\usepackage{epsfig}
\usepackage{subfigure}

\theoremstyle{plain}
\newtheorem{theorem}{Theorem}
\newtheorem{lemma}[theorem]{Lemma}

\newtheorem{proposition}[theorem]{Proposition}

\theoremstyle{definition}
\newtheorem{definition}[theorem]{Definition}

\newtheorem{problem}[theorem]{Problem}

\theoremstyle{remark}

\title{Characterizing circle graphs with binomial partial Petrial polynomials}

\author{Ruiqing Feng,
Qi Yan\footnote{Corresponding author.}, Xuan Zheng\\
\small School of Mathematics and Statistics, Lanzhou University, PR China\\
\small{\tt Email: fengrq2024@lzu.edu.cn; yanq@lzu.edu.cn; zhengx2024@lzu.edu.cn}}
\date{}
\journal{arXiv}

\begin{document}
\begin{abstract}
The partial Petrial polynomial was first introduced by Gross, Mansour, and Tucker as a generating function that enumerates the Euler genera of all possible partial Petrials on a ribbon graph. Yan and Li later extended this polynomial invariant to circle graphs by utilizing the correspondence between circle graphs and bouquets. Their explicit computation demonstrated that paths produce binomial polynomials, specifically those containing exactly two non-zero terms. This discovery led them to pose a fundamental characterization problem: identify all connected circle graphs whose partial Petrial polynomial is binomial. In this paper, we solve this open problem in terms of local complementation and prove that for connected circle graphs, the binomial property holds precisely when the graph is a path.
\end{abstract}
\maketitle

\section{Introduction}
In 1979, Wilson \cite{SEW} defined the Petrial operation for embedded graphs. This construction retains the vertex and edge sets of the original graph but replaces its faces with Petrie polygons—closed walks formed by alternating left-right traversals of edges. The Petrial admits an intuitive interpretation within the ribbon graph framework: for a ribbon graph, it is generated by detaching one end of every edge from its incident vertex disc, applying a half-twist to the edge, and reattaching the endpoint. Applying this operation to edge subsets yields the partial Petrial.
\begin{definition}[\cite{EM1}]
Let $G$ be a ribbon graph with edge set $E(G)$, and let $A \subseteq E(G)$. The \emph{partial Petrial} of $G$ with respect to $A$, denoted $G^{\times|A}$, is the ribbon graph obtained from $G$ by adding a half-twist to each edge in $A$.
\end{definition}

The growing importance of partial Petrials stems from their utility across multiple disciplines, including topological graph theory, knot theory, matroid / delta-matroid theory, and physics \cite{EM1, EM}. Gross, Mansour, and Tucker \cite{GMT2} defined the partial Petrial polynomial for arbitrary ribbon graphs, establishing formulas and recursions for families such as ladder ribbon graphs.
\begin{definition}[\cite{GMT2}]
The \emph{partial Petrial polynomial} $^{\partial}\varepsilon^{\times}_{G}(z)$ of a ribbon graph $G$ is a generating function of the numbers of partial Petrials of $G$ of the given Euler genus:
$$^{\partial}\varepsilon^{\times}_{G}(z):=\sum_{A\subseteq E(G)}z^{\varepsilon(G^{\times|A})}$$
where $\varepsilon(G^{\times|A})$ stands for the Euler genus of $G^{\times|A}$.
\end{definition}

Yan and Jin \cite{QYJ4} investigated analogues of partial Petrial polynomials for delta-matroids, proving that signed intersection graphs uniquely determine these polynomials for bouquets (single-vertex ribbon graphs). Recently, Yan and Li \cite{YQLY} established that a bouquet's partial Petrial polynomial is fundamentally determined by its intersection graph. Specifically, if two bouquets $B_1$ and $B_2$ share isomorphic intersection graphs, then
\[{^\partial\varepsilon^{\times}_{B_1}(z)}={^\partial\varepsilon^{\times}_{B_2}(z)}.\]
Furthermore, they extended the partial Petrial polynomial to circle graphs. Here, a graph is a \emph{circle graph} precisely when it is the intersection graph of a bouquet.
\begin{definition}[\cite{YQLY}]
The \emph{partial Petrial polynomial}, denoted $P_{G}^{\times}(z)$, of a circle graph $G$ is defined as
\[
P_{G}^{\times}(z) := {}^\partial\varepsilon^{\times}_{B}(z),
\]
where $B$ is a bouquet whose intersection graph is isomorphic to $G$.
\end{definition}

Yan and Li \cite{YQLY} established that for a connected graph with $n \geq 2$ vertices,
the partial Petrial polynomial has non-zero coefficients for all degrees from 1 to $n$
if and only if the graph is complete. Their explicit computation for paths yields:

\begin{theorem}[\cite{YQLY}]\label{the01}
Let $P_n$ be a path on $n \geq 1$ vertices. Then
\begin{equation*}
P_{P_n}^{\times}(z) =
\begin{cases}
\left(\frac{2^{n}-1}{3}\right)z^{n-1} + \left(\frac{2^{n+1}+1}{3}\right)z^n, & \text{if } n \text{ is even,}\\
\left(\frac{2^{n}+1}{3}\right)z^{n-1} + \left(\frac{2^{n+1}-1}{3}\right)z^n, & \text{if } n \text{ is odd.}
\end{cases}
\end{equation*}
\end{theorem}

A polynomial is called \emph{binomial} if it contains exactly two non-zero terms, i.e., is of the form $az^k + bz^m$ for $k \neq m$. They posed the characterization problem:

\begin{problem}[\cite{YQLY}]
Is a connected circle graph $G$ necessarily a path when $P_G^{\times}(z)$ is binomial? If not, characterize all such graphs.
\end{problem}

We solve this problem using local complementation, a fundamental graph operation first introduced by Kotzig \cite{KA}. This operation has significant applications across diverse domains including network science \cite{DHW, HPE, KDOS}. Our main result is:

\begin{theorem}\label{the03}
A connected circle graph $G$ on $n$ vertices has a binomial partial Petrial polynomial if and only if $G$ is a path.
\end{theorem}

\section{Preliminaries}

In this paper, we represent embedded graphs as ribbon graphs, adopting the following formal definition:

\begin{definition}[\cite{BR}]
A {\it ribbon graph} $G$ is a (orientable or non-orientable) surface with boundary,
represented as the union of two sets of topological discs, a set $V(G)$ of vertices, and a set $E(G)$ of edges, subject to the following restrictions.
\begin{enumerate}
\item[(1)] The vertices and edges intersect in disjoint line segments.
\item[(2)] Each such line segment lies on the boundary of precisely one vertex and precisely one edge.
\item[(3)] Every edge contains exactly two such line segments.
\end{enumerate}
\end{definition}

A ribbon graph is \emph{orientable} if its underlying surface is orientable; otherwise, it is \emph{non-orientable}.
An edge in a ribbon graph is called a \emph{loop} if both its endpoints are incident to the same vertex. A loop is said to be \emph{non-orientable} if the corresponding ribbon (considered together with its associated vertex) forms a M\"{o}bius band; otherwise, it is an \emph{orientable loop}.

A \emph{bouquet} is a ribbon graph with exactly one vertex. We say that two loops in a bouquet are \emph{interlaced} if their ends are met in an alternating order when travelling round the vertex boundary. The {\it intersection graph} \cite{CSL} $I(B)$ of  a bouquet $B$ is the graph with vertex set $E(B)$, and in which two vertices $e_1$ and $e_2$ of $I(B)$ are adjacent if and only if $e_1$ and $e_2$ are interlaced in $B$.

If $G$ is a ribbon graph, we denote by $f(G)$ the number of boundary components of $G$, and we define $v(G)$, $e(G)$, and $c(G)$ to be the number of vertices, edges, and connected components of $G$, respectively. We let \[\chi(G)=v(G)-e(G)+f(G),\] the usual \emph{Euler characteristic}, where $G$ is connected or not. The notation $\varepsilon(G)$ represents the \emph{Euler genus} of $G$, that is, \[\varepsilon(G)=2c(G)-\chi(G).\]

A subgraph obtained by vertex deletions only is called an \emph{induced subgraph}. If $X$ is the set of vertices deleted, the resulting subgraph is denoted by $G\setminus X$. Frequently, it is the set $Y := V\setminus X$ of vertices which remain that is the focus of interest. In
such cases, the subgraph is denoted by $G[Y]$ and referred to as the subgraph of $G$ induced by $Y$. Thus $G[Y]$ is the subgraph of $G$ whose vertex set is $Y$ and whose edge set consists of all edges of G which have both ends in $Y$.

A \emph{simple graph} is a graph without loops or multiple edges. A \emph{graft} $(G, L_G)$ consists of a simple graph $G$ and a subset $L_G\subseteq V(G)$ of vertices. The
\emph{adjacency matrix} $\mathbf{A}_{(G, L_G)}$ of a graft $(G, L_G)$ is the matrix over $GF(2)$ whose rows and columns correspond to the
vertices of $G$; and where, for $u\neq v$, the $(u, v)$-entry of $\mathbf{A}_{(G, L_G)}$ is 1 if the corresponding vertices
$u$ and $v$ are adjacent in $G$, and is 0 otherwise; and the $(v, v)$-entry of $\mathbf{A}_{(G, L_G)}$ is 1 if $v\in L_G$, and is 0 otherwise.

A fundamental graph operation is local complementation, first studied by Kotzig in \cite{KA}. We use $N_G(v)$ to denote the set of neighbours of a vertex $v$ in the graph $G$. Note that $v\notin N_G(v)$.

\begin{definition}[\cite{KA}]
Let $G$ be a simple graph and $v\in V(G)$. The \emph{local complementation} at $v$, denoted $G*v$, is the graph obtained from $G$ by replacing the induced subgraph on the neighborhood $N_G(v)$ with its complement. Equivalently, $G*v$ is formed by toggling all adjacencies between vertices in $N_G(v)$ (i.e., replacing edges with non-edges and vice versa within $N_G(v)$). We further define $G~\underline{*}~v:=(G*v)\setminus \{v\}.$
\end{definition}

The following definition adapts local complementation for grafts.

\begin{definition}[\cite{Moffatt15}]
Let $G$ be a simple graph with vertex subset $L_G \subseteq V(G)$.
For any vertex $v \in L_G$, the \emph{local complementation} at $v$ is the operation on the graft $(G, L_G)$ defined by
\[
(G, L_G) \mapsto (G \ast v,  L_G \triangle N_{G}(v)).
\]
The \emph{local complementation deletion} at $v$ is the operation
\[
(G, L_G)~\underline{*}~v := \left( G ~\underline{*}~v,  (L_G\setminus\{v\}) \triangle N_{G}(v) \right).
\]
A graft $(H, L_H)$ is a \emph{local complementation minor} of $(G, L_G)$ if it can be obtained from $(G, L_G)$ by a sequence of local complementation deletion operations.
\end{definition}

\section{Main results}
First, we establish that grafts preserve the corank of their adjacency matrices under local complementation minors.
\begin{proposition}\label{pro02}
Let $(G, L_G)$ and $(H, L_H)$ be grafts. If $(H, L_H)$ is a local complementation minor of $(G, L_G)$, then
\[
\operatorname{corank}(\mathbf{A}_{(G, L_G)}) = \operatorname{corank}(\mathbf{A}_{(H, L_H)}).
\]
\end{proposition}

\begin{proof}
Consider a local complementation deletion operation at vertex $v$ on the graft $(G, L_G)$. Define $\mathbf{A}_{(G,L_G)}^{(v)}$ as the matrix obtained from $\mathbf{A}_{(G,L_G)}$ through the following elementary operations:
\begin{itemize}
    \item Add the row corresponding to $v$ to each row indexed by vertices in $N_G(v)$ (over $GF(2)$).
    \item Add the column corresponding to $v$ to each column indexed by vertices in $N_G(v)$ (over $GF(2)$).
\end{itemize}
This transformation is shown below, where $d = |N_G(v)|$:
\[
\mathbf{A}_{(G,L_G)} =
\begin{pmatrix}
1 & \mathbf{1}^\top & \mathbf{0} \\
\mathbf{1} & \mathbf{A}_{11} & \mathbf{A}_{12} \\
\mathbf{0} & \mathbf{A}_{12}^\top & \mathbf{A}_{22}
\end{pmatrix}
\longrightarrow
\begin{pmatrix}
1 & \mathbf{1}^\top & \mathbf{0} \\
\mathbf{0} & \mathbf{A}_{11}' & \mathbf{A}_{12} \\
\mathbf{0} & \mathbf{A}_{12}^\top & \mathbf{A}_{22}
\end{pmatrix}
\longrightarrow
\begin{pmatrix}
1 & \mathbf{0}^\top & \mathbf{0} \\
\mathbf{0} & \mathbf{A}_{11}' & \mathbf{A}_{12} \\
\mathbf{0} & \mathbf{A}_{12}^\top & \mathbf{A}_{22}
\end{pmatrix}
= \mathbf{A}_{(G,L_G)}^{(v)},
\]
where $\mathbf{1}$ is the all-ones vector of length $d$, $\mathbf{A}_{11}' = \mathbf{A}_{11} + \mathbf{J}$ over $GF(2)$ ($\mathbf{J}$ is the $d \times d$ all-ones matrix), and all entries are in $\{0,1\}$.

Since elementary operations preserve matrix rank, we have
\[
\operatorname{corank}(\mathbf{A}_{(G,L_G)}^{(v)}) = \operatorname{corank}(\mathbf{A}_{(G, L_G)}).
\]

The adjacency matrix of $(G, L_G) \mathbin{\underline{*}} v$ is
\[
\mathbf{A}_{(G, L_G) \mathbin{\underline{*}} v} =
\begin{pmatrix}
\mathbf{A}_{11}' & \mathbf{A}_{12} \\
\mathbf{A}_{12}^\top & \mathbf{A}_{22}
\end{pmatrix}.
\]
From the block structure of $\mathbf{A}_{(G,L_G)}^{(v)}$, we observe
\[
\mathbf{A}_{(G,L_G)}^{(v)} =
\begin{pmatrix}
1 & \mathbf{0}^\top \\
\mathbf{0} & \mathbf{A}_{(G, L_G) \mathbin{\underline{*}} v}
\end{pmatrix}.
\]
This yields the rank relationship
\[
\operatorname{rank}(\mathbf{A}_{(G,L_G)}^{(v)}) = \operatorname{rank}(\mathbf{A}_{(G, L_G) \mathbin{\underline{*}} v}) + 1,
\]
and consequently
\[
\operatorname{corank}(\mathbf{A}_{(G, L_G) \mathbin{\underline{*}} v}) = \operatorname{corank}(\mathbf{A}_{(G,L_G)}^{(v)}).
\]
Combining these results gives
\[
\operatorname{corank}(\mathbf{A}_{(G, L_G) \mathbin{\underline{*}} v}) = \operatorname{corank}(\mathbf{A}_{(G, L_G)}).
\]

As $(H, L_H)$ is obtained from $(G, L_G)$ through a sequence of local complementation deletion operations, the equality extends inductively
\[
\operatorname{corank}(\mathbf{A}_{(G, L_G)}) = \operatorname{corank}(\mathbf{A}_{(H, L_H)}). \qedhere
\]
\end{proof}

\begin{lemma}\label{pro01}

For any path $P$ on $n$ vertices, there exists $L_P \subseteq V(P)$
such that the graft $(P, L_P)$ has a local complementation minor $(P', L_{P'})$
where $P'$ consists of a single isolated vertex and $L_{P'} = \emptyset$.
\end{lemma}

\begin{proof}
For $n = 1$, let $L_P = \emptyset$. Then $(P', L_{P'}):=(P, L_P)$ is a local complementation minor of $(P, L_P)$,
$P'$ is an isolated vertex, and $L_{P'} = L_P = \emptyset$.

For $n \geq 2$, let $V(P) = \{v_1, \dots, v_n\}$ where $v_1$ and $v_n$ are the leaves.
Set $L_P = \{v_1, v_n\}$. Performing sequential local complementation deletions at $v_1, v_2, \dots, v_{n-1}$ yields:
\[
(P', L_{P'}):=((P, L_P) ~\underline{*}~ v_1) ~\underline{*}~ v_2 ~\underline{*}~ \cdots ~\underline{*}~ v_{n-1} = (P\setminus\{v_1, \dots, v_{n-1}\}, \emptyset).
\]
Note that $P'=P\setminus\{v_1, \dots, v_{n-1}\} = \{v_n\}$ is an isolated vertex and $L_{P'} = \emptyset$.
\end{proof}

\begin{theorem}\label{the02}
Let $G$ be a simple graph that is not a path. Then there exists $L_G \subseteq V(G)$ such that the graft $(G, L_G)$ contains a local complementation minor $(G', L_{G'})$ where $G'$ has at least two isolated vertices, neither of which belongs to $L_{G'}$.
\end{theorem}

\begin{proof}
We proceed by induction on the number of vertices $n \geq 2$ of $G$. Let $V(G) = \{v_1, \dots, v_n\}$.
		
\textbf{Base case:} If $n=2$, then the only non-path graph is two isolated vertices. Set $L_G = \emptyset$. Thus $(G', L_{G'}) := (G, L_G)$ is a local complementation minor of $(G, L_G)$ with two isolated vertices not in $L_{G'} = \emptyset$.

If $n=3$, then non-path graphs are three configurations:
\begin{itemize}
    \item \textit{3-cycle}: Set $L_G = \{v_1, v_2, v_3\}$. Then
    $$(G, L_G) ~\underline{*}~ v_1 = (G~\underline{*}~v_1, \emptyset)$$
    where $G~\underline{*}~v_1$ has two isolated vertices and $L_{G'} = \emptyset$.

    \item \textit{Path $P_2$ with vertices $v_1, v_2$ union isolated vertex $v_3$}: Set $L_G = \{v_1, v_2\}$. Then
    $$(G, L_G) ~\underline{*}~ v_1 = (G~\underline{*}~v_1, \emptyset)$$
    where $G~\underline{*}~v_1$ has two isolated vertices and $L_{G'} = \emptyset$.

    \item \textit{Three isolated vertices}: Set $L_G = \emptyset$. The graft is its own minor with three isolated vertices not in $L_{G'} = \emptyset$.
\end{itemize}
In all cases, the minor has at least two isolated vertices not in $L_{G'}$.

\textbf{Inductive hypothesis:}
Assume for any non-path simple graph $H$ with $k$ vertices ($2 \leq k \leq n-1$), there exists $L_H \subseteq V(H)$ such that $(H, L_H)$ has a local complementation minor with at least two isolated vertices not in $L_H$.

\textbf{Inductive step:} Consider a non-path simple graph $G$ with $n \geq 4$ vertices.

If $G$ is disconnected, denote its connected components as $G = G_1 \cup \cdots \cup G_{c(G)}$ where $c(G) \geq 2$. For the first two components $G_j$ ($j = 1, 2$):
\begin{itemize}
    \item If $G_j$ is a path, by Lemma \ref{pro01}, there exists $L_{G_j} \subseteq V(G_j)$ such that the graft $(G_j, L_{G_j})$ has a local complementation minor $(G'_j, L_{G'_j})$
    where $G'_j$ consists of a single isolated vertex and $L_{G'_j} = \emptyset$.

    \item If $G_j$ is not a path, since $|V(G_j)| \leq n-1$, by the inductive hypothesis there exists $L_{G_j} \subseteq V(G_j)$ such that $(G_j, L_{G_j})$ has a local complementation minor $(G'_j, L_{G'_j})$ with at least two isolated vertices, neither of which belongs to $L_{G'_j}$.
\end{itemize}
Set $L_G = L_{G_1} \cup L_{G_2}$. Then the graft $\left( \bigcup_{i=1}^c G'_i, \bigcup_{i=1}^cL_{G'_i}\right)$ is a local complementation minor of $(G, L_G)$, where for $i \geq 3$, $G'_i = G_i$ and $L_{G'_i}=\emptyset.$   Moreover, $G'_1 \cup G'_2$ contains at least two isolated vertices not in $L_{G'_1} \cup L_{G'_2}$.

If $G$ is a connected non-path simple graph, it can be divided into the following two cases:

\textbf{Case 1:} $G$ is a tree.

A non-path tree must contain a branching vertex (i.e., a vertex of degree at least 3) $v \in V(G)$. Since $G$ is a tree, $G[N_G(v)]$ is an independent set. Performing local complementation deletion at $v$ yields $G~\underline{*}~v$, which is non-path (as complementing the independent set $N_G(v)$ creates a clique containing a cycle) and has $n-1$ vertices. By the inductive hypothesis, there exists $L_1 \subseteq V(G~\underline{*}~v)$ such that $(G~\underline{*}~v, L_1)$ contains a local complementation minor $(G', L_{G'})$ with at least two isolated vertices not in $L_{G'}$. Set $L_G = (L_1\cup \{v\}) \Delta N_G(v)$. Then $(G', L_{G'})$ is a local complementation minor of $(G, L_G)$.

\textbf{Case 2:} $G$ is non-tree (i.e., contains cycles).

When $G$ is a cycle $C_n$ ($n \geq 4$), select any vertex $v \in V(G)$. After local complementation deletion at $v$, $G~\underline{*}~v$ is a cycle with $n-1$ vertices. By the inductive hypothesis, there exists $L_2 \subseteq V(G~\underline{*}~v)$ such that $(G~\underline{*}~v, L_2)$ has a local complementation minor $(G', L_{G'})$ where $G'$ contains at least two isolated vertices not in $L_{G'}$. Set $L_G = (L_2 \cup \{v\}) \Delta N_G(v)$. Then $(G', L_{G'})$ is a local complementation minor of $(G, L_G)$.

When $G$ is non-tree and not a cycle (i.e., contains cycles and branching vertices), consider any vertex $v \in V(G)$:
\begin{itemize}
    \item If $G~\underline{*}~v$ is non-path, apply the inductive hypothesis directly to establish the theorem.

    \item If $G~\underline{*}~v$ is a path, retain vertex $v$ and partition $V(G)\setminus \{v\}$ into $N_v := N_G(v)$ and $N_v^c := V(G) \setminus (N_G(v)\cup \{v\})$ (note $N_v^c$ may be empty). By the definition of local complementation, both the edges between $N_v$ and $N_v^c$ and the induced subgraph $G[N_v^c]$ remain unchanged. Since $G~\underline{*}~v$ is a path, each connected component of $G[N_v^c]$ is a path, and all vertices in $N_v^c$ have degree at most 2.
    \begin{description}
  \item [Case 2.1:] If $N_v^c = \emptyset$, then $v$ is adjacent to all other vertices, i.e., $N_v = V(G)\setminus\{v\}$, implying $|N_v| = n-1 \geq 3$. Since $G~\underline{*}~v$ is a path, select a vertex $v_1 \in V(G~\underline{*}~v)$ of degree 2 (which exists since $G~\underline{*}~v$ is a path with $n-1 \geq 3$ vertices), and denote its neighbors by $x$ and $y$.
Note that $x$ and $y$ are adjacent in $G$ (as they are non-adjacent in $G~\underline{*}~v$), $v$ is adjacent to both $x$ and $y$ (since $x, y \in N_v$), and neither $x$ nor $y$ is adjacent to $v_1$ in $G$. Thus $v$, $x$, and $y$ form a 3-cycle in $G$. This 3-cycle persists in $G~\underline{*}~v_1$, implying $G~\underline{*}~v_1$ is a non-path graph with $n-1$ vertices. We apply the inductive hypothesis.
  \item [Case 2.2:] If $N_v^c \neq \emptyset$, since $G$ contains cycles and $n \geq 4$, we must have $|N_v| \geq 2$. Otherwise, if $|N_v| = 1$, then $G~\underline{*}~ v$ would contain a cycle, contradicting the assumption that $G~\underline{*}~ v$ is a path. We consider two subcases based on the size of $N_v$:
  \begin{description}
  \item [Case 2.2.1:] $|N_v| \geq 3$. Since $|N_v| \geq 3$, $v$ is a branching vertex. Select an arbitrary vertex $v_2 \in N_v^c$. As $v \notin N_G(v_2)$, $v$ remains a branching vertex in $G~\underline{*}~ v_2$, implying $G~\underline{*}~ v_2$ is non-path. The inductive hypothesis then establishes the theorem.
  \item [Case 2.2.2:] $|N_v| = 2$. Let $N_v = \{a, b\}$. Assume $a$ and $b$ are non-adjacent in $G$. Then the sequence $a-v-b$ forms a path in $G$.
Moreover, $a$ and $b$ become adjacent in $G~\underline{*}~ v$ (due to complementation of $N_v$). If exactly one of $\{a,b\}$ is adjacent to vertices in $N_v^c$, then since $G$ contains cycles, either $G[N_v^c \cup \{a\}]$ or $G[N_v^c \cup \{b\}]$ contains a cycle $C$. This cycle $C$ persists in $G~\underline{*}~ v$, contradicting its path structure. If both $a$ and $b$ are adjacent to vertices in $N_v^c$, the path structure of $G~\underline{*}~ v$ (where $a$ and $b$ are adjacent) forces $G[N_v^c]$ to consist of two disjoint paths connected to $a$ and $b$ respectively. Combined with the path $a-v-b$, this implies $G$ is itself a path, contradicting the assumption that $G$ contains a cycle. Thus $a$ and $b$ are adjacent in $G$, forming a 3-cycle $vab$. We consider two subcases based on $|V(G~\underline{*}~v)|$:
If $|V(G~\underline{*}~v)| = 3$ (i.e., $n=4$), then $G~\underline{*}~v$ is a 3-vertex path. Hence, $G~\underline{*}~ a$ consists of a 2-vertex path union an isolated vertex (hence non-path). We apply the inductive hypothesis.
If $|V(G~\underline{*}~v)| > 3$, select a vertex $c$ such that $\{a,b,c\}$ induces a disconnected subgraph in $G~\underline{*}~ v$. This implies $|N_G(c) \cap \{v,a,b\}| \leq 1$, and consequently $G~\underline{*}~ c$ contains the 3-cycle $vab$ (non-path). We apply the inductive hypothesis.
\end{description}
\end{description}
\end{itemize}
\end{proof}

\begin{lemma}[\cite{YQLY}]\label{lem05}
Let $G$ be a connected circle graph. Then the degree of the partial Petrial polynomial \(P_G^{\times}(z)\) is $|V(G)|$.
\end{lemma}

\begin{lemma}[\cite{GMT2}]\label{lem06}
For any connected circle graph $G$, the partial Petrial polynomial \(P_G^{\times}(z)\) is an interpolating polynomial.
\end{lemma}

\begin{lemma}[\cite{MB}]\label{lem03}
Let \( I(B) \) be the interlace graph of a bouquet \( B \), and let \( S \subseteq V(I(B)) \) correspond to the non-orientable loops of \( B \). Then
\[
f(B) = \operatorname{corank}\left( \mathbf{A}_{(I(B), S)} \right) + 1.
\]
\end{lemma}

\noindent\textbf{Proof of Theorem \ref{the03}.}
Sufficiency follows from Theorem~\ref{the01}. For necessity, suppose the partial Petrial polynomial \(P_G^{\times}(z)\) of the circle graph \(G\) is binomial. By Lemmas~\ref{lem05} and \ref{lem06},
\[
P_G^{\times}(z) = a_{n-1} z^{n-1} + a_n z^n,
\]
where \(a_{n-1}, a_n \neq 0\).

Since \(G\) is a circle graph, there exists an orientable bouquet \(B\) such that \(G = I(B)\). Then
\[
{}^\partial\varepsilon^{\times}_{B}(z) = P_G^{\times}(z) = a_{n-1} z^{n-1} + a_n z^n.
\]
For any subset \(D \subseteq E(B)\), \(\varepsilon(B^{\times|D})\) equals \(n-1\) or \(n\). Since \(v(B^{\times|D}) = 1\) and \(e(B^{\times|D}) = n\), Euler's formula yields
\[
f(B^{\times|D}) =2+e(B^{\times|D})-v(B^{\times|D})-\varepsilon(B^{\times|D}),
\]
implying that \(f(B^{\times|D}) = 1\) or \(2\).

As \(B\) is orientable, the set of non-orientable loops in \(B^{\times|D}\) is precisely \(D\) (corresponding to a vertex subset of \(I(B)\) under the natural bijection). By Lemma~\ref{lem03},
\begin{equation}\label{eq01}
  \operatorname{corank}(\mathbf{A}_{(G, D)}) \leq 1.
\end{equation}

Now assume \(G\) is not a path. By Theorem~\ref{the02}, there exists a vertex subset \(L_G \subseteq V(G)\) (corresponding to an edge subset of \(B\) under the natural bijection) such that the graft \((G, L_G)\) contains a local complementation minor \((G', L_{G'})\) where \(G'\) has at least two isolated vertices not in \(L_{G'}\). Note that \(\operatorname{corank}(\mathbf{A}_{(G', L_{G'})}) \geq 2\). By Proposition~\ref{pro02},
\[
\operatorname{corank}(\mathbf{A}_{(G, L_G)}) = \operatorname{corank}(\mathbf{A}_{(G', L_{G'})}) \geq 2,
\]
contradicting \eqref{eq01} (which holds for all subsets, including  \(L_G\)). Therefore, \(G\) must be a path.

\section*{Acknowledgements}
This work is supported by NSFC (Nos. 12471326, 12101600).

\end{document}